\documentclass[reqno,a4paper,12pt]{amsart}

\usepackage[utf8]{inputenc}
\usepackage[T1]{fontenc}
\usepackage[english]{babel}
\usepackage[babel]{microtype}
\usepackage[a4paper, margin=1.1in]{geometry}
\usepackage[dvipsnames]{xcolor}
\usepackage{amsmath,amssymb,amsthm}
\usepackage{csquotes}
\usepackage{bbm}
\usepackage{dsfont}
\usepackage{lmodern}
\usepackage{fixcmex}
\usepackage{cases}
\usepackage{enumitem}
\setlist[enumerate,1]{label=(\roman*)}

\usepackage{mathrsfs}
\usepackage{mathtools}
\usepackage{commath}
\usepackage{thmtools}
\usepackage{thm-restate}

\usepackage[pdfborderstyle={/S/U/W 0},unicode]{hyperref}
\hypersetup{
    colorlinks=true,
    linkcolor=blue!85!black,
    citecolor=blue!85!black,
    urlcolor=blue!85!black,
}
\usepackage{cleveref}

\numberwithin{equation}{section}

\declaretheoremstyle[
  headfont=\normalfont\bfseries,
  bodyfont=\normalfont,
  spaceabove=6pt,
  spacebelow=6pt,
  headpunct={.},
  postheadspace=1em
]{noital}

\declaretheorem[style=plain,numberwithin=section,name=Theorem]{theorem}

\declaretheorem[style=plain,sibling=theorem,name=Proposition]{proposition}

\declaretheorem[style=plain,sibling=theorem,name=Conjecture]{conjecture}

\declaretheorem[style=plain,sibling=theorem,name=Claim]{claim}

\declaretheorem[style=noital,sibling=theorem,name=Definition]{definition}

\newcommand{\cE}{\mathcal{E}}

\newcommand{\cM}{\mathcal{M}}

\newcommand{\cP}{\mathcal{P}}

\newcommand{\PP}{\mathbb{P}}

\newcommand{\ol}[1]{\overline{#1}}
\newcommand{\PPe}{\PP_{\! F}}

\begin{document}

\title[The acyclic directed bunkbed conjecture is false]{The acyclic directed bunkbed conjecture is false}

\author[T. Przybyłowski]{Tomasz Przybyłowski}

\address{Mathematical Institute, University of Oxford, Oxford, United Kingdom}
\email{tomasz.przybylowski@kellogg.ox.ac.uk}

\begin{abstract}
    We construct a simple acyclic directed graph for which the Bunkbed Conjecture is false, thereby resolving conjectures posed by Leander and by Hollom.
\end{abstract}

\maketitle


\section{Introduction}
The Bunkbed Conjecture, introduced by Kasteleyn in 1985 \cite[Remark 5]{bergkahn2001}, postulates a certain monotonicity of connectivity probabilities in Bernoulli bond percolation on an undirected graph. Being an intuitive, but surprisingly hard to prove, statement, a lot of effort has been put in proving its correctness. Despite that, the original conjecture turned out to be false, as shown recently by Gladkov, Pak and Zimin \cite{gpz2024} who built on the work of Hollom disproving the hypergraph variant \cite{hollom2024}.

Since the main conjecture does not hold, the problem of finding broad families of graphs, or other general models, where it does hold is of great interest. 
One of the remaining promising general models considered is the setup of acyclic directed graphs put forward by Leander \cite{leander2009thesis}, and highlighted by Hollom \cite{hollom2024}. Nevertheless, in this work we prove that the Bunkbed Conjecture does not hold in this case either.

Let $G=(V,E)$ be a (directed or undirected) graph, and let $(V^-, E^-)$ and $(V^+,E^+)$ be two copies of $G$. For $v \in V$ let $v^-$, $v^+$ denote the copies of $v$ in $V^-$ and $V^+$, respectively. The \emph{bunkbed graph} $\ol{G} = (\ol{V}, \ol{E})$ of $G$ is the graph defined by $\ol{V} = V^- \sqcup V^+$ and $\ol{E} = E^- \sqcup E^+ \sqcup \widetilde{E}$, where $\widetilde{E}$ consists of bidirected edges $(v^-, v^+)$ for all $v\in V$. We will call the edges in $\widetilde{E}$ \emph{vertical}, and the edges in $E^- \sqcup E^+$ \emph{horizontal}.

The \emph{bunkbed percolation model} is a probabilistic measure on the subgraphs of the bunkbed graph. The most commonly chosen measure is the $p$-Bernoulli bond percolation~$\PP_p$ of the edges of the bunkbed graph, i.e., a random subgraph of the edges of the bunkbed graph, in which every edge is retained with probability $p$, or deleted with probability $1-p$, independently of other edges. 

The original Bunkbed Conjecture states the following:
\begin{conjecture}
    For every undirected graph $G=(V,E)$ and every $u$, $v \in V$
    \[ \PP_p( u^- \leftrightarrow v^-)  \geq \PP_p( u^- \leftrightarrow v^+),\]
    where $x \leftrightarrow y$ is the event that in the percolated bunkbed graph there is a path connecting vertices $x$ and $y$.
\end{conjecture}

This conjecture was disproved by Gladkov, Pak and Zimin \cite{gpz2024}. Their counterexample is ingeniously based on a graph found by Hollom \cite{hollom2024}, which in turn served as a counterexample to the Hypergraph Bunkbed Conjecture. The apparent intuition behind the conjecture is that vertices $u^-$ and $v^-$ are \emph{closer} to each other than $u^-$ and~$v^+$ for any conceivable notion of \emph{closeness}, and so the chance of an open path between $u^-$ and $v^-$ should be at least as large as the chance of open path between $u^-$ and~$v^+$. There are certain instances where the conjecture does hold, notably on complete graphs \cite{debuyer2018,vHL19}, for adjacent vertices in edge-transitive graphs \cite{richthammer2025}, for vertices $u$, $v$ with the same neighbourhood \cite{richthammer2025}; and also for the random cluster model with parameter $q=2$ \cite{haggstrom2003probability}; for more thorough discussion see \cite{gpz2024}.

Leander considered the Bunkbed Conjecture for directed graphs in \cite{leander2009thesis}, in which case one compares the probability of existence of directed paths in a Bernoulli percolated bunkbed graph. This directed case was disproved by Hollom, who found a directed multigraph counterexample. Hollom conjectured that one can avoid the use of multiple edges. Further, Leander considered the conjecture for acyclic directed graphs. Interestingly, Hollom's counterexample contains a lot of directed cycles and as such does not disprove Leander's acyclic variant. The acyclic graphs seem promising to satisfy the conjecture as they might be amenable to some inductive argument. Despite that, we find a simple acyclic directed graph not satisfying the conjecture, thereby resolving the conjectures by Leander and by Hollom. 

\begin{theorem} \label{thm:main}
    There is a simple acyclic directed graph with distinguished vertices $u$, $v$ such that
    \[ \PP_{1/2}( u^- \rightarrow v^-) < \PP_{1/2}( u^- \rightarrow v^+).\]
\end{theorem}

\begin{figure}[h!]
\centering
\includegraphics[scale=1.3]{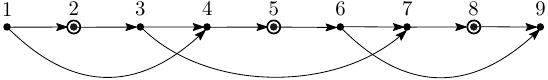}
\caption{\label{fig:main_graph}Graph $G_1$.}

\vspace{0.3cm} 

\includegraphics[scale=1.1]{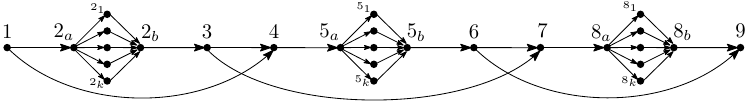}
\caption{\label{fig:gadget1}Graph $G_2^k$.}

\vspace{0.3cm} 

\includegraphics[scale=1.3]{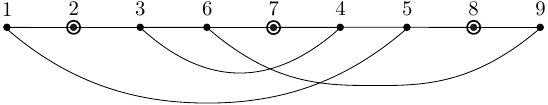}
\caption{\label{fig:hollom_graph}Hollom's counterexample to the conditioned site bunkbed conjecture.} 
\end{figure}

Our counterexample is the graph $G_2^k$ depicted in \Cref{fig:gadget1}. It is a blow-up of the graph $G_1$ depicted in \Cref{fig:main_graph}, with vertices $2$, $5$, $8$ substituted with gadgets on $k+2$ vertices. 

The proof of \Cref{thm:main} has two steps. Firstly, in \Cref{sec:cond_bond} we prove that bond percolation on $\ol{G_1}$ conditioned to have vertical edges $(2^-,2^+), (5^-,5^+), (8^-,8^+)$ does not satisfy the Bunkbed Conjecture. Then, in \Cref{sec:uncond_bond}, we show that in unconditioned bond percolation on $\ol{G_2^k}$ the gadgets behave with high probability as if they were vertical edges at $2$, $5$, $8$ in $\ol{G_1}$.

Finally, let us remark that $G_1$ bears a great resemblance to Hollom's counterexample to the conditioned site bunkbed conjecture (see \Cref{fig:hollom_graph} and also $G_2$ in Figure 2 of \cite{hollom2024}), and one can verify that the undirected version of $G_1$ itself is also a counterexample to the conditioned site bunkbed conjecture.

\section{Notation}

Let $G = (V,E)$ be a directed graph, and $\ol{G} = (\ol{V}, \ol{E}) = (V^- \sqcup V^+, E^- \sqcup E^+ \sqcup \widetilde{E})$ be its bunkbed graph, where $E^- \sqcup E^+$ are the horizontal edges in the lower and top bunk, while $\widetilde{E}$ are the vertical bidirected edges $(v^-,v^+)$ for $v \in V$. Let $\cP(\ol{E})$ denote the family of all spanning subgraphs of the bunkbed graph. This is the event space of the bunkbed percolation process, and as such we will refer to subsets of $\cP(\ol{E})$ as \emph{events}.

For a horizontal edge $e = (u^\varepsilon, v^\varepsilon) \in E^- \cup E^+$, where $\varepsilon \in \{-,+\}$, define its \emph{shadow} to be $\pi(e) = (u,v) \in E$, and let $\ol{e} = (u^{-\varepsilon}, v^{-\varepsilon})$ be its \emph{mirrored edge}. 

For $H \subseteq \ol{E}$ and $L \subseteq E$ define the \emph{$L$-mirrored subgraph} $M(H,L) \subseteq \ol{E}$ by 
\begin{equation*}
    M(H,L) = (\widetilde{E} \cap H) \cup \{e \in H \setminus \widetilde{E}\colon \pi(e) \notin L\} \cup \{\bar{e} \colon e \in H \setminus \widetilde{E} \text{ and } \pi(e) \in L\},
\end{equation*}
i.e. $M(H,L)$ consists of all vertical edges of $H$, all edges of $H$ whose shadow is not in~$L$, and the mirrored edges of those edges of $H$ whose shadow is in $L$. For a family of subgraphs $A \subseteq \cP(\ol{E})$ we let $\cM(A,L) = \{M(H,L)\colon H \in A\}$. Observe that $\cM(\cM(A,L),L) = A$, and so in particular $|A|=|\cM(A,L)|$.

We call $v \in V$ a \emph{post} in $H \subseteq \ol{E}$ if $(v^-,v^+) \in H$ and define $P(H) \subseteq V$ to be the set of all posts of $H$. 

We use the notation $\{u^\varepsilon \to v^\delta\}$ to signify the family of subgraphs of $\ol{E}$ which contain a directed path from $u^\varepsilon \in \ol{V}$ to $v^\delta \in \ol{V}$. We extend the notation and write $\{u^\pm \to v^\delta\}$ to denote the set $\{u^- \to v^\delta\} \cup \{u^+ \to v^\delta\}$, and we write $\{x,y \to z\}$ to mean $\{x \to z\} \cap \{y \to z\}$, etc. Symbols as $\{x \to y,z\}$ or $\{x,y \leftarrow z\}$ are defined similarly. The crossed arrow works as the complement, e.g. $\{u^\pm \not\to v^+\} = \cP(\ol{E}) \setminus \{u^\pm \to v^+\}$.


\section{Conditioned directed bunkbed bond model} \label{sec:cond_bond}

\begin{definition}[Conditioned]
    For a directed graph $G = (V,E)$ and a subset $T \subseteq V$, the \emph{conditioned bunkbed percolation model} $\cE_T$ is a random subgraph of $\ol{G}$ sampled uniformly from $\cP(\ol{E})$ subject to having posts at vertices $T$, i.e. $P(\cE_T) \supseteq T$.
\end{definition}

\begin{proposition} \label{prop:conditional}
    The Bunkbed Conjecture does not hold for the conditioned bunkbed model $\cE_T$ of $G_1$ with $T = \{2,5,8\}$, and $u=1$, $v=9$.
\end{proposition}

 The main strategy of proof of \Cref{prop:conditional} is to cut out from the events $1^-\to9^-$ and $1^-\to9^+$ as many equiprobable parts as possible, and then analyze the remaining refined connection events. The cut-out parts come from bijections (or \emph{mirroring}) of the respective connection events. 
 This fairly natural approach has been used many times before in the study of the problem, see e.g. \cite{linusson2011percolation,linusson2019erratum,leander2009thesis,hollom2024}. 
 
 In particular, one of the cut-out parts is the event $F^c$ of presence of an additional post apart from $2,5,8$, so that the set of posts is a cut-set between $1$ and $9$ in the undirected $G_1$. In the case of undirected graphs, if the posts form a cut-set between $u$ and $v$, then $\PP(u^- \leftrightarrow v^-) = \PP(u^- \leftrightarrow v^+)$. Interestingly, in the setup of directed graphs, it is not enough to have a \emph{directed} cut-set between $u$ and $v$ made of posts (see e.g. \Cref{prop:conditional}). Nevertheless, if the posts form a cut-set in the underlying undirected graph, then one has equality, see \cite{leander2009thesis,hollom2024}. The proof of this fact uses a fairly standard mirroring argument, which we repeat to keep the work self-contained.

\begin{proof}[Proof of \Cref{prop:conditional}]
    Throughout the proof $\PP$ is the law of the conditioned model~$\cE_T$. Let $A = \{1^- \to 9^-\}$ and $B = \{1^-\to 9^+\}$. Our goal is to show that $\PP(A) < \PP(B)$. Define $F = \{ P(\cE_T) = \{2,5,8\}\}$ and
    \begin{align*}
        A_1 &= A \cap F^c, \\
        A_2 &= A \cap F \cap \{1^- \not\to 5^\pm\}, \\
        A_3 &= A \cap F \cap \{1^- \to 5^\pm\} \cap \{1^- \not\to 3^\pm\}, \\
        A_4 &= A \cap F \cap \{1^- \to 3^-,3^+,5^\pm\},
    \end{align*}
    and define events $B_1$, \ldots, $B_4$ by substituting above $B$ in the place of $A$. Let $\widetilde\cE_T$ be the subgraph of $\cE_T$ induced by vertices $3^\pm, 5^\pm, 6^\pm, 7^\pm$, $8^\pm$, $9^\pm$. Further, define
    \begin{align*}
        P_\varepsilon &= \{1^- \to 3^\varepsilon,5^\pm\} \cap \{1^- \not\to 3^{-\varepsilon}\}, \\
        Q_{\varepsilon, \delta} &= \{ 3^\varepsilon \to 9^\delta \text{ in } \widetilde{\cE}_T\} \cup \{5^\pm \to 9^\delta\},
    \end{align*}
    for $\varepsilon, \delta \in \{-,+\}$.

    \begin{figure}[h!]
        \centering
        \includegraphics[scale=1.3]{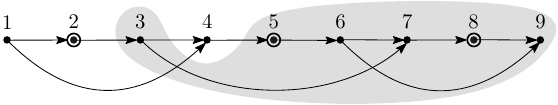}
        \caption{\label{fig:shaded}Graph $G_1$ with its shaded part corresponding to $\widetilde{\cE}_T$.}
    \end{figure}

    Observe that
    \begin{align}
        A &= A_1 \sqcup A_2 \sqcup A_3 \sqcup A_4 \sqcup (P_- \cap Q_{--} \cap F) \sqcup (P_+ \cap Q_{+-} \cap F), \nonumber \\
        B &= B_1 \sqcup B_2 \sqcup B_3 \sqcup B_4 \sqcup (P_- \cap Q_{-+} \cap F) \sqcup (P_+ \cap Q_{++} \cap F). \label{eqAB}
    \end{align}

    Let $\PPe(\cdot) = \PP(\cdot \mid F)$. We claim the following.
    \begin{claim} \label{claim:cond}\phantom{blah}
        \begin{enumerate}
            \item For $i=1,2,3,4$ we have $\PP(A_i) = \PP(B_i)$;
            \item $\PPe(P_+) > \PPe(P_-)$;
            \item $\PPe(Q_{++}) = \PPe(Q_{--}) > \PPe(Q_{-+}) = \PPe(Q_{+-})$.
        \end{enumerate}
    \end{claim}
    
    \begin{proof}[Proof of claim.]
    \textit{(i):} 
    Observe that in the underlying undirected graph of $G_1$ there is only one path connecting vertices $1$ and $9$ and avoiding $2$, $5$, $8$. Moreover, this path passes through all of the remaining vertices. Therefore, if $\cE_T$ contains a post other than $2$, $5$, $8$, then every (undirected) path connecting $1^-$ to $9^-$ or $9^+$ will necessarily pass through some post. 

    For $S \subseteq V \setminus T$ let $F_S$ be the event that $P(\cE_T) = S \cup T$. Further, let $R(S)$ be the set of edges in the underlying undirected graph of $G_1$ reachable from $9$ without passing through a vertex in $S \cup T$ (e.g., $R(\{4\}) = E \setminus \{(1,2), (1,4),(4,5)\}$, $R(\{3,9\}) = \emptyset$). 
    
    Due to the observation above, for every non-empty $S \subseteq V\setminus T$ we have \[\cM(A \cap F_S, R(S)) = B \cap F_S,\] and so $\PP(A \cap F_S) = \PP(B \cap F_S)$. Summing over all non-empty $S$ yields $\PP(A_1) = \PP(B_1)$. 
    
    Let $K = \{(5,6), (6,7), (6,9), (7,8), (8,9)\}$. We have \[\cM(A_2, \{(8,9)\}) = B_2, \quad \cM(A_3,K) = B_3 \quad\text{ and }\quad \cM(A_4,K \cup \{(3,7)\}) = B_4,\]
    and so the remaining equalities of (i) follow.

    \textit{(ii):} Let $R = \{(1^-,4^-) \in \cE_T\}$. We have 
    \[\cM(P_+ \cap R^c \cap F, \{(2,3), (3,4), (4,5)\}) = P_- \cap R^c \cap F,\]
    and so it remains to show that $\PPe(P_+ \cap R) > \PPe(P_- \cap R)$, which is true as the consideration of cases yields $\PPe(P_+ \cap R) = \tfrac{1}{16} \cdot \tfrac{5}{8}$ and $\PPe(P_- \cap R) = \frac{1}{16}\cdot\tfrac{1}{2}$.

    \textit{(iii)}: The two equalities are true due to the mirroring of all of the edges. Further,
    \begin{align}
        \PPe(Q_{--}) = \PPe(3^- \to 9^- \text{ in } \widetilde\cE_T) + \PPe(5^\pm \to 9^-) - \tfrac{1}{2} \PPe(7^-,5^\pm \to 9^-), \label{eq:qmm}\\
        \PPe(Q_{-+}) = \PPe(3^- \to 9^+ \text{ in } \widetilde\cE_T) + \PPe(5^\pm \to 9^+) - \tfrac{1}{2} \PPe(7^-,5^\pm \to 9^+), \label{eq:qmp}
    \end{align}
    where we use the equality
    \begin{align*}
    \PPe(3^- \to 9^- \text{ in } \widetilde\cE_T \text{ and } 5^\pm \to 9^-) &= \PPe((3^-,7^-) \in \cE_T \text{ and } 7^-,5^\pm\to 9^-) \\
    &= \tfrac{1}{2}\PPe(7^-,5^\pm\to 9^-),
    \end{align*}
    and an analogous one for $\PPe(3^- \to 9^+ \text{ in } \widetilde\cE_T \text{ and } 5^\pm \to 9^+)$.

    The mirroring of $(8,9)$, and of all of the edges yields that the first two terms in \eqref{eq:qmm} and \eqref{eq:qmp} are equal. Thus, it is enough to show $\PPe(7^-,5^\pm \to 9^-) < \PPe(7^+,5^\pm \to 9^-)$, which is equivalent to
    \begin{equation} \label{eq:preiso}
        \PPe(7^-,5^\pm \to 9^- \text{ and } 7^+ \not\to 9^-) < \PPe(7^+,5^\pm \to 9^- \text{ and } 7^- \not\to 9^-).
    \end{equation}
    Observe that reversing all arrows in $G_1$ and mapping each vertex $i \mapsto 10-i$ gives an isomorphic copy of $G_1$. Further, if a subgraph $H$ satisfies $x^\varepsilon \to y^\delta$, then the corresponding subgraph after mapping satisfies $(10-y)^\delta \to (10-x)^\varepsilon$, etc. Hence, under this map, the event $\{7^-,5^\pm \to 9^- \text{ and } 7^+ \not\to 9^-\}$ becomes $\{3^-,5^\pm \leftarrow 1^- \text{ and } 3^+\not\leftarrow 1^-\}$, while $\{7^+,5^\pm \to 9^- \text{ and } 7^- \not\to 9^-\}$ becomes $\{3^+,5^\pm \leftarrow 1^- \text{ and } 3^-\not\leftarrow 1^-\}$. Thus, \eqref{eq:preiso} is equivalent to $\PPe(P_-) < \PPe(P_+)$ from claim~(ii), which completes the proof of (iii).\qedhere%
    \end{proof}

    Finally, from \eqref{eqAB}, the claim, and the independence of $P_\varepsilon \mid F$ and $Q_{\varepsilon, \delta} \mid F$, we get
    \begin{align*}
        \PP(&B) - \PP(A) \\
        &= \PP(F) \left( \PPe( P_+) \PPe(Q_{++}) + \PPe(P_-)\PPe(Q_{-+}) - \PPe(P_-)\PPe(Q_{--}) - \PPe(P_+)\PPe(Q_{+-}) \right) \\
        &= \PP(F) (\PPe(P_+) - \PPe(P_-))(\PPe(Q_{--}) - \PPe(Q_{-+})) > 0,
    \end{align*}
    as desired.
\end{proof}


\section{Unconditioned directed bunkbed bond model} \label{sec:uncond_bond}

    \begin{definition}
        For a directed graph $G = (V,E)$, the \emph{(unconditioned) bunkbed percolation model} $\cE$ is a random subgraph of $\ol{G}$ sampled uniformly from $\cP(\ol{E})$.
    \end{definition}

We boost the conditioned counterexample to an unconditioned one, thereby proving \Cref{thm:main}. In order to do so, we substitute in the graph $G_1$ each of the vertices $i = 2,5,8$ with a gadget on $k+2$ vertices $i_a,i_b,i_1,\ldots,i_k$, as shown in \Cref{fig:gadget1}. We call the resulting graph $G_2^k$. 


\begin{proof}[Proof of \Cref{thm:main}]
Let $A_k = \{1^- \to 9^-\} \subseteq \cP(\ol{E(G^k_2)})$ and $B_k = \{1^- \to 9^+\} \subseteq \cP(\ol{E(G^k_2)})$. We will show that for sufficiently large $k$ the percolation model $\cE$ satisfies
\begin{equation}
    \PP(A_k) < \PP(B_k). \label{eq:G2k}
\end{equation}

For $i = 2,5,8$ and $j \in [k]$ let $C^i_j$ be the event that all the five edges $(i_j^-, i_j^+)$, $(i_a^\varepsilon, i_j^\varepsilon)$, $(i_j^\varepsilon, i_b^\varepsilon)$ for $\varepsilon \in \{-,+\}$ are present in $\cE$. Further, let $C_i = \bigcup_j C^i_j$, and $C = C_2 \cap C_5 \cap C_8$. We have $\PP(C_i) = 1 - (1-2^{-5})^k$ and so $\PP(C) \to 1$ as $k \to \infty$.

For every subgraph $H \subseteq \ol{E(G^k_2)}$ in $C$ define $\gamma(H) \subseteq \ol{E(G_1)}$ by contracting in~$H$ all the vertices $i_a^\varepsilon$, $i_b^\varepsilon$, $i_1^\varepsilon$, \ldots, $i_k^\varepsilon$ to a single vertex $i^\varepsilon$ for all $i=2,5,8$ and $\varepsilon \in \{-,+\}$ (and erasing potential loops). Observe that if $H$ is chosen uniformly from $C$, then $\gamma(H)$ is distributed as $\cE_T$ in \Cref{prop:conditional}. Further, for $H \in C$ the event $1^-\to9^\varepsilon$ holds for~$H$ if and only if $1^-\to9^\varepsilon$ holds for $\gamma(H)$.

As a result,
\begin{align*}
\PP(A_k) &= \PP(C) \cdot \PP(A_k \mid C) + \PP(C^c \cap A_k) \\
&= \PP(C) \cdot \PP(1^- \to 9^- \text{ in } \cE_T) + \PP(C^c \cap A_k) \stackrel{k\to\infty}{\longrightarrow} \PP(1^- \to 9^- \text{ in } \cE_T) = \PP(A),
\end{align*}
and similarly $\PP(B_k) \to \PP(B)$ as $k$ tends to $\infty$. Thus, due to \Cref{prop:conditional}, for sufficiently large~$k$ inequality \eqref{eq:G2k} holds, as desired.
\end{proof}

\section*{Acknowledgements}
   I would like to thank Oliver Riordan for comments on a previous draft. This work was supported by the Additional Funding Programme for Mathematical Sciences, delivered by EPSRC (EP/V521917/1) and the Heilbronn Institute for Mathematical Research.


\bibliographystyle{amsplain}
\providecommand{\bysame}{\leavevmode\hbox to3em{\hrulefill}\thinspace}
\providecommand{\MR}{\relax\ifhmode\unskip\space\fi MR }
\providecommand{\MRhref}[2]{%
  \href{http://www.ams.org/mathscinet-getitem?mr=#1}{#2}
}
\providecommand{\href}[2]{#2}


\end{document}